\documentclass[12pt]{article}
\newcommand{\hf}{\textstyle{\frac{1}{2}}}

\usepackage{natbib,amsmath,amssymb,amsthm,amsfonts}
\usepackage{graphicx,url}
\usepackage{latexsym}

\newcommand{\mR}{\mathbb{R}}
\newcommand{\mC}{\mathbb{C}}

\theoremstyle{plain}

\newtheorem{proposition}{Proposition}

\begin{document}

\title{Joint density of eigenvalues in spiked multivariate models}

\author{Prathapasinghe Dharmawansa, Iain
  M. Johnstone \\ 
  Department of Statistics, Stanford U.} 




\maketitle

\begin{abstract}
  The classical methods of multivariate analysis are based on the
  eigenvalues of one or two sample covariance matrices.  In many
  applications of these methods, for example to high dimensional data,
  it is natural to consider alternative hypotheses which are a low
  rank departure from the null hypothesis.  For rank one alternatives,
  this note provides a representation for the joint eigenvalue density
  in terms of a single contour integral. This will be of use for
  deriving approximate distributions for likelihood ratios and
  `linear' statistics used in testing.
\end{abstract}



\section{Introduction}

The eigenvalues of one or two sample covariance matrices play a
central role in multivariate analysis. A long list of
examples, including principal components analysis (PCA), canonical
correlations analysis (CCA), multivariate analysis of variance
(MANOVA) and multiple response linear regression are the main subject
of many standard textbooks, such as \citet{mkb79,ande03}.

Under the common assumption of Gaussian data, much is known about the
joint and marginal distribution of the eigenvalues.
For example, under the typical null hypotheses, the joint density of the
eigenvalues has an explicit formula, derived in 1939 in the celebrated
and independent work of
Fisher, Girshick, Hsu, Mood and Roy.
Under general alternatives, the joint density is given by an integral
over a group of matrices. If the number of variables, and hence
eigenvalues, is large, $p$ say, as is common nowadays, this integral
will be high dimensional, of dimension $O(p^2)$. 

A remarkable classification of the joint density functions was given
by \citet{jame64}, using hypergeometric functions of matrix argument.
He showed how the classical multivariate methods could be organized
into five cases, involving hypergeometric functions $\,_pF_q$ of
different orders, specifically 
$\,_0F_0, \,_0F_1, \,_1F_0, \,_1F_1,$ and $\,_2F_1.$ 
Remarkable though this work is, and despite significant progress on
the numerical computation of hypergeometric functions, e.g. 
\citet{koed06}, these expressions for the joint densities have proved
challenging to work with in application.

In many high dimensional applications, however, it may be reasonable
to consider alternative hypotheses which are low rank departures from
the null. For some examples, see \citet{jona14}.
In this note we consider the simplest case, namely rank one deviations, and
show that the joint eigenvalue density can then be reduced to a single
(contour) integral.

We believe this integral representation to be of interest at least
because it is amenable to approximation when dimension $p$ is large,
leading to simple approximations to at least certain aspects of these
multivariate eigenvalue distributions.

We mention two examples of such applications.
\begin{itemize}
\item[(i)] derivation of limiting Gaussian approximations for `linear
  statistics' (including, for example, the likelihood ratio test, and
  `high-dimension-corrected' likelihood ratio test,
  \citet{omh13,wsy13}). Particular
  cases ($\,_0F_0, \,_0F_1, \,_1F_1$) have been given for complex data
  by \citet{pmc14}.
\item[(ii)] delineation of the region of contiguous alternatives to
  the null hypothesis, and description of the Gaussian limit for the
  log-likelihood ratio process inside the contiguity region.
This leads to a comparative understanding of the power properties of
various hypothesis tests, both traditional and new, in the contiguity
region. This example has been studied in the case of PCA,
corresponding to $\,_0F_0$, by \citet{omh13}, and
work is in progress to apply the result of this note to the general 
$\,_pF_q$ cases.
\end{itemize}

We will adopt James' systematization in order to give a unified
derivation of our contour formulas.
We give the rank one formula for $\,_pF_q$ in real and complex cases,
Section 2.
This can be converted directly into an expression for the joint
density function for the eigenvalues in each of James' five cases
(for both $\mR$ and $\mC$).  Section 3 illustrates this process in
one case, testing equality of covariance matrices, for real data 
(i.e. $\,_1F_0$).

In the real case, the proof of Section \ref{sec:cont-integr-repr}
applies only to even dimension $p$. Section
\ref{sec:real-case-integer} gives a different proof valid for all
integer $p$.

\section{Contour integral representation for rank one}
\label{sec:cont-integr-repr}

Let $X, Y$ be $r \times r$ Hermitian matrices.
The definitions of hypergeometric functions with one and two matrix
arguments are given, for example, by \citet{jame64}, with separate
expressions for real and complex cases.

The definitions simplify in our special case in which $X$ has rank
one, with nonzero eigenvalue $x$. 
For $a \in \mC$, let $(a)_k = a(a+1) \cdots (a+k-1), (a)_0 = 1$ be the rising
factorial, and for vectors of parameters
$a = (a_l)_{l=1}^p, b = (b_l)_{l=1}^q$
with $a_l \in \mC$ and $b_l \in \mC \backslash \{0, -1, -2, \ldots \}$, 
adopt the abbreviation
\begin{displaymath}
  \rho_k = \rho_k(a,b) 
    = \frac{ (a_1)_k \cdots (a_p)_k}{(b_1)_k \cdots (b_q)_k}.
\end{displaymath}
If $X$ has rank one as described, define
\begin{equation}
  \label{eq:rankone}
  \,_pF_q^\alpha (a,b; X, Y)
  = \sum_{k=0}^\infty \rho_k \frac{(1/\alpha)_k}{(r/\alpha)_k}
         \frac{x^k C_k^\alpha(Y)}{k!}.
\end{equation}
Here $\alpha > 0$ indexes a one parameter family that includes the
real ($\alpha = 2$) and complex ($\alpha = 1$) cases. 
Also, $C_k^\alpha$ are Jack polynomials (e.g. \citet{macd95}): in the 
real case ($\alpha = 2$), they
reduce to James' zonal polynomials (e.g. \citet{muir82}), and in the complex
case ($\alpha=1$), to a normalization of the Schur functions
(e.g. \citet{des07}). 
A contour formula for $C_k^\alpha (Y)$ is quoted below; for now we
note that $C_k^\alpha (X) = x^k$, and (e.g. \citet[eq. (245)]{wang12}) that 
\begin{displaymath}
  C_k^\alpha(I) = \prod_{j=0}^{k-1} \frac{r+\alpha j}{1+\alpha j}
     = \frac{(r/\alpha)_k}{(1/\alpha)_k},
\end{displaymath}
which explains the form of the two ratios in formula
(\ref{eq:rankone}) as
$C_k^\alpha(X) C_k^\alpha(Y) / C_k^\alpha(I)$.

The series (\ref{eq:rankone}) converges for all $x, Y$ if $p \leq q$; 
for $x \| Y \| < 1$ if $p=q +1$ (and $\|Y\|$ denotes the maximum
eigenvalue in absolute value of $Y$) and finally diverges unless
it terminates if $p > q+1$
(e.g. \citet{mph95}).

With this notation, the scalar generalized hypergeometric function,
which does not depend on $\alpha$, is
\begin{displaymath}
  \,_pF_q (a,b;x) = \sum_{k=0}^\infty \rho_k(a,b) \frac{x^k}{k!}.
\end{displaymath}

The main result of this note can now be stated.
\begin{proposition}
Suppose that $p \leq q+1$, 
$X$ is rank 1 with positive eigenvalue $x$ and that $Y$
is positive definite with eigenvalues $(y_j)_{j=1}^r$.

(i) Suppose that $r/\alpha$ is a positive integer, say $r/\alpha =
m+1$, and that  
$a_l \notin \{1, \ldots, m\}$ and $ b_l \notin \{ m, m-1, m-2, \ldots \}$.
Then, 
\begin{equation}
  \label{eq:contour-formula}
    \,_pF_q^\alpha (a,b; X, Y)
    = \frac{\Gamma(m+1)}{x^m \rho_m'} \frac{1}{2 \pi i} 
    \int_\Gamma \,_pF_q (a-m,b-m; xs) 
     \prod_{j=1}^r \frac{1}{(s-y_j)^{1/\alpha}}  {\rm d}s,
\end{equation}
where the contour $\Gamma$ starts at $- \infty$, encircles $0$ and $\{
y_j \}$ counterclockwise and returns to $-\infty$. Further, 
$a-m$ denotes the vector with entries $a_i -m$ and
\begin{displaymath}
  \rho_m' = \rho_m(a-m, b-m).
\end{displaymath}
Equality holds in the common domain of analyticity of both sides: 
$\mC$ if $p \leq q$ and $\mC \backslash (1,\infty)$ if $p = q+1$. 

(ii) If instead $r/\alpha = m + \epsilon$ for $\epsilon \in (0,1)$ and
non-negative integer $m$, then under the same conditions
\begin{equation}
  \label{eq:contour-formula-mod}
    \,_pF_q^\alpha (a,b; X, Y)
    = \frac{(\epsilon)_m}{x^m \rho_m'} \frac{1}{2 \pi i} 
    \int_\Gamma s^{\epsilon -1} \,_{p+1}F_{q+1} (a-m,1,b-m,\epsilon; xs) 
     \prod_{j=1}^r \frac{1}{(s-y_j)^{1/\alpha}}  {\rm d}s.
\end{equation}
(iii) If $\alpha = 2$, then formula (2) holds for \textbf{any} integer
$r$, still with $m+1 = r/2$, if 
the symbol $(a)_m$ is interpreted as $\Gamma(a+m)/\Gamma(a)$ for
non-integer $m$.
\end{proposition}

Thus, in the real ($\alpha = 2$) and complex ($\alpha = 1$) cases of most
interest in applications, formula (\ref{eq:contour-formula}) holds for
all positive integer $r$. 

Particular cases of (\ref{eq:contour-formula}) are already known:
$\,_0F_0$ for both real and complex cases 
\citep{mo12,omh13}, for general $\alpha$, \citet{wang12,forr11},
and for the complex
case only, $\,_0F_1$ \citep{pd13}
and $\,_1F_1$ \citep{pmc14}.
\citet{wang12} also gives formula (\ref{eq:contour-formula-mod}) in the 
$\,_0F_0$ case.
A generalization of (i)
to the multi-spike case
has been given for $\,_0F_0$ by \citet{onat14} and recently extended
to $\,_pF_q$ by \citet{pmc14a}.

\textbf{Proof.}
  Parts (i) and (ii) are shown here; part (iii) uses a different
  argument and is deferred to Section 4.
  We begin with a result from \citet[eq. (248)]{wang12}, which states that
\begin{displaymath}
(1/\alpha)_k  \, \frac{C_k^\alpha(Y)}{k!}
  = \frac{1}{2 \pi i} \int_{\Gamma'} \prod_{j=1}^r
  \frac{1}{(1-zy_j)^{1/\alpha}}  \frac{{\rm d}z}{z^{k+1}}.    
\end{displaymath}
Here the contour $\Gamma'$ encircles zero and is chosen small enough so that all
$y_j^{-1}$ lie outside.

Insert this into \eqref{eq:rankone} and interchange summation and
integration to obtain 
\begin{equation}
  \label{eq:Fcontour}
  \,_pF_q^\alpha (a,b; X, Y)
   = \frac{1}{2 \pi i} \int_{\Gamma'} \prod_{j=1}^r
   \frac{1}{(1-zy_j)^{1/\alpha}} G(z;x) {\rm d}z 
\end{equation}
where the series
\begin{displaymath}
  G(z;x) = \sum_{k=0}^\infty \frac{\rho_k}{(r/\alpha)_k}
  \frac{x^k}{z^{k+1}}
\end{displaymath}
converges for all $x, z$ if $p \leq q$ and for $|x/z| < 1$ if $p =
q+1$.

Now write $r/\alpha = m+1$ and introduce the variable $l = k+m$, so
that
\begin{equation}
  \label{eq:secondform}
  G(z;x) = \sum_{l=m}^\infty \frac{\rho_{l-m}}{(m+1)_{l-m}}
               \frac{x^{l-m}}{z^{l-m+1}}
         = \frac{m!}{x^m} \frac{z^{m-1}}{\rho_m'}
           \sum_{l=m}^\infty \frac{\rho_l(a-m,b-m)}{l!}
           \left(\frac{x}{z}\right)^l, 
\end{equation}
where we have used $(m+1)_{l-m} = l!/m!$, and noted 
that $(c)_{l-m} = (c-m)_l / (c-m)_m$ so that
\begin{displaymath}
  \rho_{l-m}(a,b) = \frac{ \rho_l(a-m,b-m)}{\rho_m(a-m,b-m)}.
\end{displaymath}
Let $G_0(z;x)$ denote the function obtained by extending the summation
in \eqref{eq:secondform} down to $l = 0$, so that
\begin{displaymath}
  G_0(z;x) = \frac{m!}{x^m} \frac{z^{m-1}}{\rho_m'} 
             \,_pF_q (a-m,b-m; x/z).
\end{displaymath}
Since we are adding a
polynomial to $G$ and a term that is analytic within the contour in
\eqref{eq:secondform}, the value of the integral is unchanged.
Hence
\begin{displaymath}
  \,_pF_q^\alpha (a,b; X, Y)
  = \frac{m!}{x^m} \frac{1}{\rho_m'}  
  \frac{1}{2 \pi i} \int_{\Gamma'} \prod_{j=1}^r
   \frac{z^{m-1}}{(1-zy_j)^{1/\alpha}} \,_pF_q (a-m,b-m; x/z) {\rm d}z.
\end{displaymath}

The change of variables $z = 1/s$ yields
\begin{align*}
  \frac{1}{2 \pi i} \int_{\Gamma'} \frac{z^{m-1}}{\prod (1-zy_j)^{1/\alpha}}
  F(x/z) {\rm d}z 
  & = \frac{1}{2 \pi i} \int_{\Gamma''} \frac{1}{s^{m+1}} 
\frac{F(xs)}{\prod (1-y_j/s)^{1/\alpha}} {\rm d}s, \\
  & = \frac{1}{2 \pi i} \int_\Gamma \frac{F(xs)}{\prod (s-y_j)^{1/\alpha}}
   {\rm d}s, 
\end{align*}
where the image $\Gamma''$ of $\Gamma'$ is deformed to $\Gamma$ as
described in the Proposition statement in order to avoid the branch
cut in the final formula. Here we use the analytic continuations of
$\,_pF_q$: entire for $p \leq q$ and for $p = q+1$ analytic off the
positive real axis $(1,\infty)$. 
The result follows.

When $r/\alpha = m + \epsilon$, we modify the argument.
In (\ref{eq:secondform}), replace $(m+1)_{l-m}$ by 
$(m+\epsilon)_{l-m}=(\epsilon)_l/(\epsilon)_m$ to obtain
\begin{displaymath}
    G(z;x) = \frac{(\epsilon)_m}{x^m} \frac{z^{m-1}}{\rho_m'}
           \sum_{l=m}^\infty \frac{\rho_l(a-m,b-m)(1)_l}{(\epsilon)_l}
           \frac{1}{l!} \left(\frac{x}{z}\right)^l.
\end{displaymath}
Proceeding as before, and extending the summation to $l=0$, so that 
\begin{displaymath}
  G_0(z;x) = \frac{(\epsilon)_m}{x^m} \frac{z^{m-1}}{\rho_m'} 
             \,_{p+1}F_{q+1} (a-m,1,b-m,\epsilon; x/z),
\end{displaymath}
we obtain formula (\ref{eq:contour-formula-mod}).

\section{Example}
\label{sec:example}

Consider the problem of testing equality of covariance matrices---the
$\,_1F_0$ case in \citet{jame64}.
Thus, suppose that $n_1, n_2 \geq p$ and that
$p \times n_1$ and $p \times n_2$ real data matrices 
$X = [ X_1 \cdots X_{n_1}]$ and 
$Y = [ Y_1 \cdots Y_{n_2}]$ have columns $X_\nu, Y_\nu$ with mean zero and
covariance matrices $\Sigma_1$ and $\Sigma_2$ respectively.
A signal detection application is described in
\citet[Sec. 3]{jona14}. 

Suppose that the observation vectors are independent Gaussian, so that 
$A_1 = X X'$ and $A_2 = Y Y'$ have Wishart distributions 
$W_p(n_1,\Sigma_1)$ and $W_p(n_2,\Sigma_2)$ respectively.
Then \citet[eq. (65)]{jame64} gives an expression for the joint
density of the eigenvalues $(f_j)$ of $A_1 A_2^{-1}$. 
To state it, we introduce notation 
$|A| = \det(A), F = \text{diag} (f_j)$ and $\Delta = \Sigma_1
\Sigma_2^{-1}$. 
We transform this expression, following \citet[p. 313-4]{muir82}, to
obtain for $n=n_1+n_2$ and  $f_1 > f_2 > \cdots > f_p$, 
\begin{equation}
  \label{eq:1F0-eq}
  p(f;\Delta)
   = \frac{c_{p,n_1,n_2}}{|\Delta|^{n_1/2}} 
     \frac{|F|^{(n_1-p-1)/2}}{|I+F|^{n/2}}
     \,_1F_0(\tfrac{n}{2};I-\Delta^{-1}, F(I+F)^{-1})
     \prod_{j<j'}^p (f_j - f_{j'}),
\end{equation}
where in this real case, $\alpha = 2$, we have written
$\,_1F_0$ for $\,_1F_0^2$.
The normalization constant is given in terms of the multivariate gamma
function \citep[p. 61]{muir82} by
\begin{displaymath}
  c_{p,n_1,n_2}
  = \frac{\pi^{p^2/2} \Gamma_p(\hf n)}{\Gamma_p(\hf p) \Gamma_p(\hf
    n_1) \Gamma_p(\hf n_2)}.
\end{displaymath}

In the spirit of application (ii) in the Introduction, we may consider the
likelihood ratio for testing the null hypothesis that $\Sigma_1 =
\Sigma_2$. Writing $\Lambda = F(I+F)^{-1}$, we have
\begin{displaymath}
  L(\Delta; \Lambda)
  = \frac{p(\Lambda; \Delta)}{p(\Lambda; I)}
  = |\Delta|^{-n_1/2} \,_1F_0(\tfrac{n}{2}; I- \Delta^{-1}, \Lambda).
\end{displaymath}

Turning now to apply the result of this paper, suppose that $\Sigma_1$
is a rank one perturbation of $\Sigma_2$, so that 
$\Sigma_1 = (I + \psi h \psi') \Sigma_2$ for real $h$ and for $\psi$ a
unit vector in $\mathbb{R}^p$. 
In this case, $\Delta = I + \psi h \psi'$, so that 
$I - \Delta^{-1}$ has rank one, with nonzero eigenvalue
$\tau = h/(1+h)$. 

Since all components of $\Lambda = F(I+F)^{-1}$ are less than one, we
may apply the contour formula (\ref{eq:contour-formula}).
Since $\,_1F_0 (a;x) = (1-x)^{-a}$, we obtain
\begin{displaymath}
  L(\tau; \Lambda)
     = \frac{n-p}{2} B\biggl(\frac{p}{2}, \frac{n-p}{2}\biggr) 
      \frac{(1-\tau)^{n_1 /2}}{\tau^{p/2 -1}} 
    \frac{1}{2 \pi i} \int_\Gamma \frac{(1-\tau s)^{-(n-p+2)/2}}{\prod_j
    (s-\lambda_j)^{1/2}} {\rm d}s.
\end{displaymath}
where $B(\alpha, \beta)$ is the usual beta function.
This is a form suitable for asymptotic approximation, the details of
which will be reported elsewhere.

\textit{Remark.} \ 
A useful check on this last formula is obtained by letting the error
degrees of freedom $n_2 \to \infty$ while keeping $p$ and $n_1$ fixed.
This limit corresponds to the case where $\Sigma_2$ is known, say
$\Sigma_2 = I$ for convenience here, and we consider the single matrix
rank one model $\Sigma_1 = I + \psi h \psi'$ and test the hypothesis
that $h = 0$. 
To compare with the formula of \citet[Lemma 3]{omh13}, let
$(\mu_j)$ be the eigenvalues of $n_1^{-1} A_1 (n_2^{-1} A_2)^{-1}$, so
that $\lambda_j = f_j/(1+f_j) = n_1 \mu_j/(n_2 + n_1 \mu_j)$.
With the change of variables $s = n_1 z/n_2$,
the previous display converges to
\begin{displaymath}
  L(\tau; \mu)
    = \Gamma\left(\frac{p}{2}\right) \left(\frac{2}{n_1}\right)^{p/2-1}
      \frac{(1-\tau)^{n_1 /2}}{\tau^{p/2 -1}} 
      \frac{1}{2 \pi i} \int_{\Gamma} e^{n_1 \tau z/2} 
       \prod_j (z - \mu_j)^{-1/2} {\rm d}z,
\end{displaymath}
which is the cited expression for the $\,_0F_0$ likelihood ratio.


\section{Real Case, integer $r$}
\label{sec:real-case-integer}

Here we prove Proposition $1$, for {\it real} matrices with integer
dimension $r$, \textit{not necessarily even}.
A similar result, with proof extending that of \citet[Lemma 2]{omh13}
has been obtained by Alexei Onatski (personal communication)  and will
appear elsewhere. 

Our goal is to prove the validity of the following
expression for $0 \leq p \leq q+1$:
\begin{align}
 \label{master1}
 {}_pF^2_q(a,b;X,Y)=\frac{\Gamma(m+1)}{x^m\rho_m'}\frac{1}{2\pi i}\int_\Gamma
 {}_pF_q\left(a-m,b-m;xs\right) \Delta_y(s) {\rm d}s
 \end{align}
where we have defined $\Delta_y(s) =
\prod_{j=1}^r\left(s-y_j\right)^{-\frac{1}{2}}$. 
The contour $\Gamma$ starts from $-\infty$ and encircles
$y_1,y_2, \cdots,y_r$ in the positive direction (i.e.,
counter-clockwise) and goes back to $-\infty$. 

In what follows, we provide an inductive proof for the above claim. 
First we establish the initial cases: ${}_0F_q$ for $q \geq 0$ and,
separately, ${}_1F_0$.
The inductive step establishes truth for ${}_{p+1}F_{q+1}$ given truth
for ${}_pF_q$. 
Also, it is worth mentioning that we assume all powers have
their principal values and all angles in the range $[-\pi,\pi)$. 

The following alternative representation of the hypergeometric
function of two matrix arguments is useful in the sequel. Let
$\mathcal{O}(r)$ be the orthogonal group and let $({\rm d} Q)$ be the
invariant measure on $\mathcal{O}(r)$ normalized to make the total
measure unity. Then, following \cite{jame64}, we can write
\begin{align}
 \label{mxint}
 {}_pF^2_{q}\left(a,b;X,Y\right)=\int_{\mathcal{O}(r)} {}_pF^2_{q}\left(a, b;X Q'Y Q\right) ({\rm d} Q).
 \end{align}
 Moreover, let us assume, without loss of generality, that $Y=\text{diag}\left(y_1,y_2,\cdots,y_r\right)$.
Since $X$ is rank-$1$, we can further simplify (\ref{mxint}) to yield
\begin{align}
\label{matrix_int}
 {}_pF^2_{q}\left(a,b;X,Y\right)
= \int_{\mathcal{S}(r)} {}_pF_{q}\left(a,b;x q_r'Y q_r\right) ({\rm d}q_r).
 \end{align}
where $\mathcal{S}(r)$ is the $r-1$ dimensional sphere embedded in $\mathbb{R}^r$,  $q_r$ is the first column of $Q$ and $({\rm d}q_r)$ is the invariant measure on $\mathcal{S}(r)$ normalized such that the total measure is one.

\subsection{Initial cases}
We first show that the statement
 (\ref{master1}) is true for ${}_0F_q$. 
With the standard notation 
${}_0\mathbf{F}_q(b;z) = {}_0F_q(b;z)/\prod_{j=1}^q \Gamma(b_j)$,
this is equivalent to showing that, for $q\geq 0$, 
\begin{equation}
 \label{master3}
 {}_0\mathbf{F}^2_{q}(b;X,Y)=\frac{\Gamma(m+1)}{x^m} \frac{1}{2\pi i}\int_{\Gamma}
 {}_0\mathbf{F}_{q}\left(b-m;xs\right) \Delta_y(s) {\rm d}s.
 \end{equation}
Our tool is a contour representation of
\cite[eq. (7.4)]{erde37}:
 \begin{align}
 \label{erd_eq}
 {}_0\mathbf{F}_q(b;z)=\frac{1}{(2\pi i)^q}
 \int_{-\infty}^{(0+)}\cdots \int_{-\infty}^{(0+)} e^{\left(\sum_{j=1}^q w_j+\frac{z}{\prod_{j=1}^q w_j}\right)}\prod_{j=1}^q \frac{{\rm d}w_j}{w_j^{b_j}}
 \end{align}
where each contour starts from $-\infty$ and encircles the origin in
the positive sense and goes back to $-\infty$. 
We use multi-index notation $w^b = \prod w_j^{b_j}, w = \prod w_j$ and
${\rm d}w = \prod {\rm d}w_j$.

We use the spherical average (\ref{matrix_int}), then Erd\'elyi's
representation, 
and change order of integration, to get
\begin{align}
 \label{beforeint}
  {}_0\mathbf{F}^2_q(b;X,Y)
   & =  \int_{\mathcal{S}(r)} {}_0\mathbf{F}_{q}\left(b;xq_r'Y
     q_r\right) ({\rm d}q_r) \\
   & = \frac{1}{(2\pi i)^q} 
  \int_{-\infty}^{(0+)}\cdots \int_{-\infty}^{(0+)}
  e^{\sum_{j=1}^q w_j}
  \int_{\mathcal{S}(r)} 
  e^{\frac{x}{w} q_r'Y q_r }
   ({\rm d}q_r) \frac{{\rm d}w}{w^b}. 
 \end{align}
A change of variable in \cite[Lemma 2]{omh13} shows that for $x,w > 0$,
\begin{equation}
  \label{eq:onatski}
   \int_{\mathcal{S}(r)} 
  e^{\frac{x}{w} q_r'Y q_r} ({\rm d}q_r)
  = \frac{\Gamma(r/2)}{ 2 \pi i} \Bigl( \frac{w}{x} \Bigr)^{r/2 -1} 
     \int_\Gamma e^{\frac{x}{w}s}\Delta_y(s) 
      {\rm d}s, 
\end{equation}
and the equality extends by analyticity to all nonzero $w \in \mC$.
Inserting this integral in (\ref{beforeint}) and noting that $\frac{r}{2}=m+1$, we  obtain
 \begin{align*}
{}_0\mathbf{F}^2_q(b;X,Y)=& \frac{\Gamma(m+1)}{x^m (2\pi i)^{q+1}}
\int_{-\infty}^{(0+)}\cdots \int_{-\infty}^{(0+)}
e^{\sum_{j=1}^q w_j}
\int_{\Gamma} 
e^{\frac{xs}{w}} \Delta_y(s)
{\rm d}s 
\frac{{\rm d}w}{w^{b-m}}
 \end{align*}
Finally, we change the order of integration and again make use of
(\ref{erd_eq}) to arrive at the desired equality \eqref{master3}.
This proves the validity of the statement (\ref{master1}) for $p=0$. 

\bigskip
\bigskip

Now we show that, for
$x\max\{y_j\}<1$, 
\begin{align}
\label{1f0main}
{}_{1}F^2_{0}(a;X,Y) =\frac{\Gamma(m+1)}{x^m (a-m)_m}\frac{1}{2\pi i}
\int_{\Gamma}  {}_{1}F_{0}\left(a-m;xs\right) \Delta_y(s)  
{\rm d}s.
\end{align}

We use identity (\ref{matrix_int}), the special form
${}_1F_0(a;z)=(1-z)^{-a}$ and the relation
\begin{equation}
\label{eq:gamma}
\frac{1}{s^a}=\frac{1}{\Gamma(a)}\int_0^\infty t^{a-1} e^{-st} {\rm
  d}t, \qquad \;\; \Re(s)>0, \Re(a)>0
\end{equation}
to obtain, after observing that $x\max\{y_j\}<1$ implies $x q_r' Y
q_r<1$, 
\begin{equation}
\label{eq_beg}
{}_1F^2_0\left(a;X,Y\right)
  = \int_{\mathcal{S}(r)}\frac{1}{\left(1-x q_r' Y q_r\right)^a}({\rm d} q_r)
  = \int_0^\infty t^{a-1} e^{-t} 
    \int_{\mathcal{S}(r)} e^{txq_r' 
       Y q_r} ({\rm d} q_r) \;{\rm d}t.
\end{equation}
Now substitute the contour identity (\ref{eq:onatski}) with $t = 1/w$,
and with the contour chosen to encircle $\{ y_j \}$ and to lie to the
left of $1/x$. We obtain
\begin{align*}
  {}_1F^2_0\left(a;X,Y\right)
   & =\frac{\Gamma\left(r/2\right)}{\Gamma(a)
      {x}^{\frac{r}{2}-1}} \frac{1}{2\pi i} 
      \int_0^\infty \int_{\Gamma}  t^{a-\frac{r}{2}} e^{-t\left(1-x
          s\right)} \Delta_y(s) {\rm d}s\; {\rm d}t \\
   & =\frac{\Gamma\left(r/2\right)\Gamma\left(a+1-\frac{r}{2}\right)}{\Gamma(a)
     x^{\frac{r}{2}-1}} 
       \frac{1}{2\pi i}
       \int_{\Gamma} (1-xs)^{r/2 -a-1}   
\Delta_y(s)  {\rm d}s,
\end{align*}
valid for $\Re(a)>\frac{r}{2}-1$, after changing order of integration
and using (\ref{eq:gamma}) and the fact that $\Re(s) < 1/x$. 
Recalling that $m = r/2 - 1$ and ${}_1F_0(a;z) = (1-z)^{-a}$, the
final form reduces to the right hand side of (\ref{1f0main}), 
under the condition $\Re(a)>m$. However, the both sides of the above
equality,  which we have established only in the domain $\Re(a)>m$ of
complex plane, are analytic functions. Therefore, the equality must
hold in the whole region of the analyticity of $a$. This establishes
the claim (\ref{1f0main}).

\subsection{Inductive step}
\label{sec:inductive-step}

First, some notation. 
We write $a_+ = (\alpha, a_1, \ldots, a_p)$ and
$b_+ = (\beta, b_1, \ldots, b_q)$ for the augmentations of $a$ and
$b$, and abbreviate ${}_{p+1}F_{q+1}$ by ${}_{p+}F_{q+}$.
Thus, the induction step 
amounts to establishing the validity of the following statement, given
the statement (\ref{master1}) is true 
\begin{equation}
 \label{ini}
 {}_{p+}F^2_{q+}(a_+,b_+;X,Y) 
  =\frac{\Gamma(m+1)}{x^m \rho_{m+}'}\frac{1}{2\pi i}\int_{\Gamma}
 {}_{p+}F_{q+}\left(a_+ -m, b_+ -m;xs\right) \Delta_y(s) {\rm d}s
 \end{equation}
 where 
 \begin{equation}
\label{eq:rhomplus}
\rho_{m+}'= \rho_m'
\frac{\Gamma(\alpha)\Gamma(\beta-m)}{\Gamma(\alpha-m) \Gamma(\beta)}.
 \end{equation}
We use a reparametrized version of the beta density
\begin{equation*}
  \phi(t; \alpha,\beta)
   = \frac{\Gamma(\beta)}{\Gamma(\alpha) \Gamma(\beta-\alpha)}
   t^{\alpha-1} (1-t)^{\beta-\alpha-1},
\end{equation*}
and the
integral representation of
the generalized hypergeometric function \citep[eq. (3.2)]{erde37} 
 \begin{equation}
 \label{hypo_alt}
 {}_{p+}F_{q+}(a_+,b_+;x)
  = \int_{0}^1 \phi(t; \alpha, \beta) \,
   {}_{p}F_{q}(a,b;xt) {\rm d}t
 \end{equation}
 where $\Re(\beta)>\Re(\alpha)>0$,
along with (\ref{matrix_int}),
 in order to write the left side of
 (\ref{ini}) as 
 \begin{align*}
{}_{p+}F^2_{q+}(a_+,b_+;X,Y) 
& = \int_{\mathcal{S}(r)} 
{}_{p+}F_{q+}(a_+,b_+;x \, q_r'Y q_r) ({\rm d} q_r)\nonumber\\
& = \int_{\mathcal{S}(r)} \int_0^1 \phi(t; \alpha, \beta)
 {}_pF_{q}(a,b;xt \, q_r'Y q_r){\rm d}t \;({\rm d} q_r) \nonumber\\
& = \int_0^1  \phi(t; \alpha, \beta) \int_{\mathcal{S}(r)} 
 {}_pF_{q}(a,b;xt \, q_r'Y q_r)({\rm d} q_r) \;{\rm d}t \\
& = \int_0^1  \phi(t; \alpha, \beta) \;  {}_pF^2_{q}\left(a,b;tX,Y\right)
\;{\rm d}t,
 \end{align*}
where we have changed the order of integration and again used
(\ref{matrix_int}). 
The final expression can be rewritten with the help of
our induction hypothesis (\ref{master1}) as
 \begin{equation}
\label{eq:rep}
 \frac{\Gamma(m+1)}{x^m\rho_m'}\frac{1}{2\pi i}
 \int_0^1 t^{-m} \phi(t; \alpha, \beta)
  \int_{\Gamma}
 {}_pF_{q}\left(a-m,b-m;xts\right) \Delta_y(s)  {\rm d}s\; {\rm d}t.
 \end{equation}
Now use the identity
\begin{equation*}
  t^{-m} \phi(t; \alpha, \beta)
   = \phi(t; \alpha-m, \beta-m) \frac{\Gamma(\beta)
     \Gamma(\alpha-m)}{\Gamma(\beta-m)\Gamma(\alpha)} 
\end{equation*}
and note from \eqref{eq:rhomplus} that the ratio of Gamma functions
equals $\rho_m'/\rho_{m+}'$. Inserting this into \eqref{eq:rep} and 
changing the order of integration, we obtain
 \begin{equation*}
 \frac{\Gamma(m+1)}{x^m\rho_{m+}'}\frac{1}{2\pi i}
 \int_{\Gamma} \Delta_y(s) 
 \int_0^1  \phi(t; \alpha-m,\beta-m)
 {}_pF_{q}\left(a-m,b-m;xts\right) {\rm d}t\; {\rm d}s.
 \end{equation*}
Now again use (\ref{hypo_alt}), along 
with the restriction $\Re(\alpha)>m$, to yield \eqref{ini}
in the domain $\Re(\beta)>\Re(\alpha)>m$ of $\mC$.
Since both sides of equality \eqref{ini}
are analytic functions, the equality must hold in the whole region of
the analyticity of $\alpha$ and $\beta$. 
This completes the induction step.



\textbf{Acknowledgements.}
This work was supported by the
  Simons Foundation Math + X program (PD)
and NIH grant 5R01 EB 001988.

\bibliographystyle{agsm}
\bibliography{contour}

\end{document}